\documentclass[11pt,leqno]{article}
\usepackage{amsfonts}
\usepackage{latexsym}
\usepackage{amsmath}
\usepackage{amssymb}
\usepackage{amsthm}
\newtheorem{theorem}{Theorem}[section]
\newtheorem{lemma}[theorem]{Lemma}

\newtheorem{proposition}[theorem]{Proposition}
\newtheorem{corollary}[theorem]{Corollary}
\newtheorem{remark}[theorem]{Remark}

\theoremstyle{definition}

\numberwithin{equation}{section}
\theoremstyle{remark}

\newtheorem*{note*}{Note}

\makeatletter
\newcommand{\ls}{\leqslant}
\newcommand{\gr}{\geqslant}
\makeatother

\begin{document}

\small

\title{\bf The isotropic position and the reverse Santal\'{o} inequality}

\author{A. Giannopoulos, G. Paouris and B-H. Vritsiou}

\date{}

\maketitle

\begin{abstract}
\footnotesize We present proofs of the reverse Santal\'{o}
inequality, the existence of $M$-ellipsoids and the reverse
Brunn--Minkowski inequality, using purely convex geometric tools.
Our approach is based on properties of the isotropic position.
\end{abstract}

\section{Introduction}

We work in ${\mathbb R}^n$, which is equipped with a Euclidean
structure $\langle\cdot ,\cdot\rangle $. We denote the corresponding
Euclidean norm by $\|\cdot \|_2$, and write $B_2^n$ for the
Euclidean unit ball, and $S^{n-1}$ for the unit sphere. Volume is
denoted by $|\cdot |$.

\smallskip

A convex body $K$ in ${\mathbb R}^n$ is a compact convex subset of
${\mathbb R}^n$ with non-empty interior. We say that $K$ is
symmetric if $x\in K$ implies that $-x\in K$. We say that $K$ is
centered if its barycenter is at the origin, i.e. $\int_K\langle
x,\theta\rangle \,d x=0$ for every $\theta\in S^{n-1}$. For every
interior point $x$ of $K$, we define the polar body $(K-x)^{\circ }$
of $K$ with respect to $x$ as follows:
\begin{equation}(K-x)^{\circ}:=\{y\in {\mathbb R}^n: \langle z-x,y\rangle \ls 1\;\hbox{for all}\; z\in K\}.\end{equation}
Note that $(K-x)^{\circ\circ}= K-x$.

\medskip

The purpose of this article is to present an alternative route to
some fundamental theorems of the asymptotic theory of convex bodies:
the reverse Santal\'{o} inequality, the existence of $M$-ellipsoids
and the reverse Brunn--Minkowski inequality. The starting point for
our approach is the isotropic position of a convex body, which can
be shown to simultaneously be an $M$-position for the body if its
isotropic constant is bounded. The new ingredient in this paper is a
way to also show, using only basic tools from the theory of convex
bodies and log-concave measures, that every convex body with bounded
isotropic constant satisfies the reverse Santal\'{o} inequality, and
then that all bodies do.

We first recall the statements and the history of the results. The
classical Blaschke-Santal\'{o} inequality states that for every
symmetric convex body $K$ in ${\mathbb R}^n$, the volume
product $s(K):=|K||K^{\circ }|$ is less than or equal to the volume product
$s(B_2^n)$, and equality holds if and only if $K$ is an ellipsoid. More generally,
for every convex body $K$, there exists a unique point $z$ in the
interior of $K$ such that
\begin{equation}|(K-z)^{\circ}|=\inf_{x\in {\rm int}(K)}|(K-x)^{\circ}|,\end{equation}
and for this point we have
\begin{equation}|K||(K-z)^{\circ}|\ls s(B_2^n)\end{equation}
(with equality again if and only if $K$ is an ellipsoid). This
unique point is usually called the Santal\'{o} point of $K$ and is
characterized by the following property: the polar body
$(K-z)^{\circ}$ of $K$ with respect to the point $z$ has its
barycenter at the origin if and only
if $z$ is the Santal\'{o} point of $K$. Observe now that the body
$K-{\rm bar}(K)$ is centered and it is the polar body of $(K-{\rm
bar}(K))^{\circ}$ with respect to the origin, hence $0$ is the
Santal\'{o} point of $(K-{\rm bar}(K))^{\circ}$. This means that for
every centered convex body $K$,
\begin{equation}s(K)=|K||K^{\circ}|=\inf_{x\in {\rm int}(K^{\circ})}|K^{\circ}||(K^{\circ}-x)^{\circ}|,\end{equation}
and this allows us to restate the Blaschke-Santal\'{o} inequality
in a more concise way: for every centered convex body $K$ in ${\mathbb R}^n$, $\,s(K)\ls s(B_2^n)$,
with equality if and only if $K$ is an ellipsoid.

In the opposite direction, a well-known conjecture of Mahler states
that $s(K)\gr 4^n/n!$ for every symmetric convex body $K$, and that
$s(K)\gr (n+1)^{n+1}/(n!)^2$ in the not necessarily symmetric case.
This has been verified for some classes of bodies, e.g. zonoids and
1-unconditional bodies (see \cite{Reisner}, \cite{Meyer},
\cite{SaintR} and \cite{GMR}). The reverse Santal\'{o} inequality,
or the Bourgain--Milman inequality, tells us that there exists an
absolute constant $c>0$ such that
\begin{equation}\left (\frac{s(K)}{s(B_2^n)}\right )^{1/n}\gr c\end{equation}
for every convex body $K$ in ${\mathbb R}^n$ which contains $0$ in
its interior. The inequality was first proved in \cite{BM} and
answers the question of Mahler in the asymptotic sense: for every
centered convex body $K$ in ${\mathbb R}^n$, the affine invariant
$s(K)^{1/n}$ is of the order of $1/n$. A few other proofs have
appeared (see \cite{iso-sym}, \cite{Kup}, \cite{Naz}), the most
recent of which give the best lower bounds for the constant $c$ and
exploit tools from quite diverse areas: Kuperberg in \cite{Kup}
shows that in the symmetric case we have $c\gr 1/2$, and his proof
uses tools from differential geometry, while Nazarov's proof
\cite{Naz} uses multivariable complex analysis and leads to the
bound $c\gr \pi^2/32$. It should also be mentioned that Kuperberg
had previously given an elementary proof \cite{Kup2} of the weaker
lower bound $s(K)^{1/n}\gr c/(n\log n)$.

The original proof of the reverse Santal\'{o} inequality in
\cite{BM} employed a dimension descending procedure which was based
on Milman's quotient of subspace theorem. Thus, an essential tool
was the $MM^{\ast }$-estimate which follows from Pisier's inequality
for the norm of the Rademacher projection. In \cite{iso-sym}, Milman
offered a second approach, which introduced an ``isomorphic
symmetrization" technique. This is a symmetrization scheme which is
in many ways different from the classical symmetrizations. In each
step, none of the natural parameters of the body is being preserved,
but the ones which are of interest remain under control. The
$MM^{\ast }$-estimate is again crucial for the proof.

\medskip

Our approach is based on properties of the isotropic position of a
convex body and combines a very simple one-step isomorphic
symmetrization argument (which is reminiscent of \cite{iso-sym})
with the method of convex perturbations that Klartag invented in
\cite{perturbation} for his solution to the isomorphic slicing
problem. Aside from the use of the latter, the approach is 
elementary, in the sense that it uses only standard tools from convex geometry; 
namely, some classical consequences of the Brunn--Minkowski inequality. Recall that a
convex body $K$ in ${\mathbb R}^n$ is called isotropic if it has
volume $1$, it is centered and its inertia matrix is a multiple of
the identity: there exists a constant $L_K >0$ such that
\begin{equation}\label{isotropic-condition1}
\int_K\langle x,\theta\rangle^2dx =L_K^2\end{equation} for every
$\theta\in S^{n-1}$. It is relatively easy to show that every convex
body has an isotropic position and that this position is
well-defined (by this we mean unique up to orthogonal transformations): if
$K$ is a centered convex body, then any linear image $\tilde{K}$ of
$K$ which has volume 1 and satisfies
\begin{equation}\label{isotropic-condition2}
\int_{\tilde{K}}\|x\|_2^2\,dx = \inf\Bigl\{\int_{T(\tilde{K})}\|x\|_2^2\,dx:T\ \hbox{is linear and volume-preserving}\Bigr\}
\end{equation}
is an isotropic image of $K$. This also implies that any isotropic
image of $K$ has the same isotropic constant, and thus $L_K$ can
be defined for the entire affine class of $K$. One of the main problems in the
asymptotic theory of convex bodies is the hyperplane conjecture,
which, in an equivalent formulation, says that there exists an absolute
constant $C>0$ such that
\begin{equation}L_n:=\max\{ L_K:K\ \hbox{is isotropic in}\ {\mathbb
R}^n\}\ls C.\end{equation} A classical reference on the subject is
the paper of Milman and Pajor \cite{MP} (see also \cite{Gian}). The
problem remains open: Bourgain \cite{Bourgain} has obtained the
upper bound $L_K\ls c\sqrt[4]{n}\log n$, and Klartag
\cite{perturbation} has improved that to $L_K\ls c\sqrt[4]{n}$ --
see also \cite{KEM}. However, in this paper we only need a few basic
results from the theory of isotropic convex bodies and, more
generally, of isotropic log-concave probability measures. All this
background information is given in Section 2; there we also list a
few more necessary tools from the general asymptotic theory of
convex bodies and, in order to stress the fact that all of them are
of purely ``convex geometric nature", we include a
short description of the arguments leading to them.

\smallskip

In Section 3 we prove the reverse Santal\'{o} inequality in two
stages. First, using elementary covering estimates, we prove a
version of it which involves the isotropic constant $L_K$ of $K$.

\begin{theorem}
Let $K$ be a convex body in $\mathbb R^n$ which contains $0$ in its
interior. Then
\begin{equation}4ns(K)^{1/n}\gr ns(K\!-\!K)^{1/n}\gr \frac{c_1}{L_K},
\end{equation}where $c_1>0$ is an absolute constant.
\end{theorem}

Then, we use Klartag's ideas from \cite{perturbation} to show that
every symmetric convex body $K$ is ``close" to a convex body $T$
with isotropic constant $L_T$ bounded by $1/\sqrt{ns(K)^{1/n}}$.

\begin{theorem}Let $K$ be a symmetric convex body in ${\mathbb R}^n$.
There exists a convex body $T$ in ${\mathbb R}^n$ such that {\rm
(i)} $c_2K\subseteq T\!-\!T\subseteq c_3K$ and {\rm (ii)} $L_T\ls
c_4/\sqrt{ns(K)^{1/n}}$, where $c_2, c_3, c_4 >0$ are absolute constants.
\end{theorem}

Since $K$ and $T\!-\!T$ have bounded geometric distance, we easily
check that $s(K)^{1/n}\simeq s(T\!-\!T)^{1/n}$. Then we can use
Theorem 1.1 for $T$ to obtain the lower bound $L_T\gr
c_5/\bigl(ns(K)^{1/n}\bigr)$. Combining this estimate with Theorem 1.2(ii), we
immediately get the reverse Santal\'{o} inequality for symmetric
bodies, and hence for all bodies.

\begin{theorem}
Let $K$ be a symmetric convex body in $\mathbb R^n$. Then
\begin{equation}s(K)^{1/n}\gr\frac{c_6}{n},
\end{equation}where $c_6>0$ is an absolute constant.
\end{theorem}

\smallskip

In Section 4 we briefly indicate how one can use Theorem 1.3 in
order to establish the existence of $M$-ellipsoids and the reverse
Brunn--Minkowski inequality. The procedure is rather standard.

The existence of an ``$M$-ellipsoid" associated with any centered
convex body $K$ in ${\mathbb R}^n$ was proved by Milman in
\cite{Mil4} (see also \cite{iso-sym}): there exists an absolute
constant $c>0$ such that for any centered convex body $K$ in
${\mathbb R}^n$ we can find an origin symmetric ellipsoid ${\cal
E}_K$ satisfying $|K|=|{\cal E}_K|$ and
\begin{gather}\label{Mellipsoid}
\frac{1}{c}|{\cal E}_K+T|^{1/n}\ls |K+T|^{1/n} \ls c |{\cal E}_K+T|^{1/n},
\\ \nonumber
\frac{1}{c}|{\cal E}_K^{\circ }+T|^{1/n}\ls |K^{\circ }+T|^{1/n}\ls
c|{\cal E}_K^{\circ }+T|^{1/n},\end{gather} for every convex body
$T$ in ${\mathbb R}^n$. The existence of $M$-ellipsoids can be
equivalently established by introducing the $M$-position of a convex
body. To any given centered convex body $K$ in ${\mathbb R}^n$ we
can apply a linear transformation and find a position
$\tilde{K}=u_K(K)$ of volume $|\tilde{K}|=|K|$ such that
(\ref{Mellipsoid}) is satisfied with ${\cal E}_K$ a multiple of
$B_2^n$. This is the so-called $M$-position of $K$. It follows then
that for every pair of convex bodies $K_1$ and $K_2$ in ${\mathbb
R}^n$ and for all $t_1,t_2>0$,
\begin{equation}\label{reverseBM}|t_1\tilde{K_1}+t_2\tilde{K_2}|^{1/n}\ls c'\left (
t_1|\tilde{K_1}|^{1/n}+t_2|\tilde{K_2}|^{1/n}\right ),\end{equation}
where $c'>0$ is an absolute constant, and that (\ref{reverseBM})
remains true if we replace $\tilde{K_1}$ or $\tilde{K_2}$ (or both)
by their polars. This statement is Milman's reverse Brunn-Minkowski
inequality.

Another way to define the $M$-position of a convex body
is through covering numbers. Recall that the covering number $N(A,B)$ of
a body $A$ by a second body $B$ is the least integer $N$ for which there exist $N$
translates of $B$ whose union covers $A$. Then, as Milman proved, there exists an absolute constant
$\beta >0$ such that every centered convex body $K$ in ${\mathbb R}^n$ has a linear image
$\tilde{K}$ which satisfies $|\tilde{K}|=|B_2^n|$ and
\begin{equation}\label{betaMposition}
\max\{ N(\tilde{K},B_2^n), N(B_2^n,\tilde{K}),N(\tilde{K}^{\circ
},B_2^n), N(B_2^n,\tilde{K}^{\circ })\} \ls\exp(\beta
n).\end{equation} We say that a convex body $K$ which satisfies
(\ref{betaMposition}) is in $M$-position with constant $\beta $. If
$K_1$ and $K_2$ are two such convex bodies, there is a standard way
to show that they and their polar bodies satisfy the reverse
Brunn--Minkowski inequality (\ref{reverseBM}) (see the end of
Section 4). Note that $M$-ellipsoids and the $M$-position of a
convex body are not uniquely defined; see \cite{Bobkov} for a recent
description in terms of isotropic restricted Gaussian measures.

\smallskip

Pisier (see \cite{Pisier-crelle} and \cite[Chapter 7]{Pisier}) has
proposed a different approach to these results, which allows one to
find a whole family of special $M$-ellipsoids satisfying stronger
entropy estimates. The precise statement is as follows. For every
$0<\alpha <2$ and every symmetric convex body $K$ in ${\mathbb
R}^n$, there exists a linear image $\tilde{K}$ of $K$ which
satisfies $|\tilde{K}|=|B_2^n|$ and
\begin{equation}\label{alpha-regular}\max\{ N(\tilde{K},tB_2^n),N(B_2^n,t\tilde{K}),
N(\tilde{K}^{\circ },tB_2^n), N(B_2^n,t\tilde{K}^{\circ })\}\ls\exp
\left (\frac{c(\alpha )n}{t^{\alpha }}\right )\end{equation} for
every $t\gr 1$, where $c(\alpha )$ is a constant depending only on
$\alpha $, with $c(\alpha )=O\big ((2-\alpha )^{-1}\big )$ as
$\alpha\to 2$. We then say that $\tilde{K}$ is in $M$-position of
order $\alpha $ (or $\alpha $-regular $M$-position). It is an
interesting question to give an elementary proof of the existence
of, say, an $1$-regular $M$-position. Another interesting question is
to check if the isotropic position is $\alpha $-regular for some
$\alpha\gr 1$ (assuming that $L_K\simeq 1$).

\section{Tools from asymptotic convex geometry}

\noindent {\bf 2.1. Basic notation}. As mentioned at the beginning of the Introduction,
we denote the Euclidean norm on ${\mathbb R}^n$ by $\|\cdot\|_2$. More generally,
if $K$ is a convex body in ${\mathbb R}^n$ which contains 0 in its interior, then we
write $p_K$ for its Minkowski functional which is defined as follows:
\begin{equation}p_K(x):=\inf\{r>0:x\in rK\},\ \, x\in {\mathbb R}^n.\end{equation}
If $K$ is symmetric, we also write $\|\cdot\|_K$ instead of $p_K$.
For every $q\gr 1$ and every symmetric convex body $B$, we define
\begin{equation}I_q(K,B):= \left(\frac{1}{|K|^{1+\frac{q}{n}}}\int_K \|x\|_B^q \, dx\right)^{1/q}.\end{equation}
If $B$ is the Euclidean ball $B_2^n$ and $K$ is an isotropic convex
body in ${\mathbb R}^n$, then from (\ref{isotropic-condition1}) we
see that
\begin{equation}I_2^2(K,B_2^n) = \int_K \|x\|_2^2 \, dx =
\int_K\Bigl(\sum_{i=1}^n \langle x,e_i\rangle^2\Bigr) \, dx = nL_K^2,\end{equation}
so $L_K = I_2(K,B_2^n)/\sqrt{n}$. More generally, as was explained in the Introduction,
if $K$ is an arbitrary convex body in ${\mathbb R}^n$, and we write $\tilde{K}$ for the
translate of $K$ which is centered, $\tilde{K}= K\!-\!{\rm bar}(K)$, then the isotropic constant
$L_K$ of $K$ can be defined by
\begin{equation}L_K:= \frac{1}{\sqrt{n}}\,\inf\bigl\{I_2\bigl(T(\tilde{K}), B_2^n\bigr): T \ \hbox{is an invertible linear transformation}\bigr\}.\end{equation}
In the sequel, we write $\overline{B}$ for the homothetic image of volume $1$ of a convex body
$B\subset \mathbb R^n$, i.e. $\overline{B}:=\frac{B}{|B|^{1/n}}$.

\medskip

As a generalization to convex bodies, we also consider
logarithmically concave (or log-concave) measures on ${\mathbb
R}^n$. This more general approach is justified by a well-known and
very fruitful idea of K. Ball from \cite{Ball} which allows one to
transfer results from the setting of convex bodies to the broader
setting of log-concave measures and vice versa. We write
${\mathcal{P}}_{[n]}$ for the class of all Borel probability
measures on $\mathbb R^n$ which are absolutely continuous with
respect to the Lebesgue measure. The density of $\mu \in
{\mathcal{P}}_{[n]}$ is denoted by $f_{\mu}$. A probability measure
$\mu \in {\mathcal{P}}_{[n]} $ is called symmetric if $f_{\mu}$ is
an even function on $\mathbb R^n$. We say that $\mu \in
{\mathcal{P}}_{[n]}$ is centered if  for all $\theta\in S^{n-1}$,
\begin{equation}
\int_{\mathbb R^n} \langle x, \theta \rangle d\mu(x) = \int_{\mathbb
R^n} \langle x, \theta \rangle f_{\mu}(x) dx = 0.
\end{equation}
A measure $\mu$ on $\mathbb R^n$ is called $\log$-concave if for any
Borel subsets $A$ and $B$ of ${\mathbb R}^n$ and any $\lambda \in
(0,1)$, $\mu(\lambda A+(1-\lambda)B) \gr
\mu(A)^{\lambda}\mu(B)^{1-\lambda}$. A function $f:\mathbb R^n
\rightarrow [0,\infty)$ is called $\log$-concave if $\log{f}$ is
concave on its support $\{f>0\}$. It is known that if a probability
measure $\mu $ is log-concave and $\mu (H)<1$ for every hyperplane
$H$, then $\mu \in {\mathcal{P}}_{[n]}$ and its density $f_{\mu}$ is
$\log$-concave (see \cite{Borell}). Note that if $K$ is a convex
body in $\mathbb R^n$, then the Brunn-Minkowski inequality implies
that ${\bf 1}_{K} $ is the density of a $\log$-concave measure.

There is also a way to generalize the notion of the isotropic constant
of a convex body in the setting of $\log$-concave measures. Set
\begin{equation}\|\mu\|_{\infty}=\sup_{x\in {\mathbb R}^n} f_{\mu} (x).\end{equation}
The isotropic constant of $\mu $ is defined by
\begin{equation}\label{definition-isotropic}
L_{\mu }:=\left (\frac{\|\mu\|_{\infty }}{\int_{{\mathbb R}^n}f_{\mu}(x)dx}\right )^{\frac{1}{n}}
[\det {\rm Cov}(\mu)]^{\frac{1}{2n}},\end{equation}
where ${\rm Cov}(\mu)$ is the covariance matrix of $\mu$ with entries
\begin{equation}{\rm Cov}(\mu )_{ij}:=\frac{\int_{{\mathbb R}^n}x_ix_j f_{\mu}
(x)\,dx}{\int_{{\mathbb R}^n} f_{\mu} (x)\,dx}-\frac{\int_{{\mathbb
R}^n}x_i f_{\mu} (x)\,dx}{\int_{{\mathbb R}^n} f_{\mu}
(x)\,dx}\frac{\int_{{\mathbb R}^n}x_j f_{\mu} (x)\,dx}{\int_{{\mathbb
R}^n} f_{\mu} (x)\,dx}\end{equation}
(in the case that $\mu$ is a centered probability measure, we can write more simply
${\rm Cov}(\mu )_{ij}:=\int_{{\mathbb R}^n}x_ix_j f_{\mu}(x)\,dx$). It is
straightforward to see that this definition coincides with the original
definition of the isotropic constant when $f_{\mu}$ is the characteristic
function of a convex body. In addition, any bounds that we have for the isotropic
constants of convex bodies continue to hold essentially in this more general setting.
This can be seen through the following construction:
let $\mu \in {\mathcal{P}}_{[n]}$ and assume that $0 \in {\rm supp}(\mu )$.
For every $p>0$, we define a set $K_p(\mu )$ as follows:
\begin{equation}\label{Kp-definition}K_p(\mu):= \left \{ x\in \mathbb R^n: p\int_{0}^{\infty}
f_{\mu}(rx)r^{p-1} dr \gr f_{\mu}(0) \right \}.\end{equation} The
sets $K_p(\mu )$ were introduced in \cite{Ball} and allow us to
study log-concave measures using convex bodies. K. Ball proved that
if $\mu$ is $\log$-concave, then $K_{p}(\mu)$ is a convex body.
Moreover, if $ \mu$ is centered, then $ K_{n+1}(\mu)$ is also
centered, and we can prove that
\begin{equation}c_1L_{K_{n+1}(\mu)}\ls L_{\mu} \ls c_2L_{K_{n+1}(\mu)}\end{equation}
for some constants $c_1, c_2 >0$ independent of $n$.

\medskip

For basic facts from the Brunn-Minkowski theory and the asymptotic theory of
finite dimensional normed spaces, we refer to the books \cite{Sch}, \cite{MS} and \cite{Pisier}.

The letters $c,c^{\prime }, c_1, c_2$ etc. denote absolute positive
constants whose value may change from line to line. Whenever we
write $a\simeq b$ for two quantities $a,b$ associated with convex bodies
or measures on ${\mathbb R}^n$, we mean that we can find positive constants
$c_1,c_2$, independent of the dimension $n$, such that $c_1a\ls b\ls c_2a$.
Also, if $K,L\subseteq \mathbb R^n$, we will write $K\simeq L$ if there exist
absolute positive constants $c_1, c_2$ such that $c_{1}K\subseteq L \subseteq c_{2}K$.

\bigskip

In the rest of the section, we collect several tools and results
from the asymptotic theory of convex bodies which will be used in
Section 3.

\medskip

\noindent {\bf 2.2. Some lemmas on covering numbers}.
Let $K, B$ be convex bodies in ${\mathbb R}^n$ with $B$ symmetric. We will give
an estimate for the covering numbers $N(K,tB)$, $t>0$, in terms
of the quantity
\begin{equation}I_1(K,B)=\frac{1}{|K|^{1+\frac{1}{n}}}\int_K\|x\|_B \, dx.\end{equation}

\smallskip

\begin{lemma}\label{covering-K-by-tB}Let $K$ be a convex body
of volume $1$ in $\mathbb R^n$ containing $0$ as an interior point.
For any symmetric convex body $B$ in ${\mathbb R}^n$ and any $t>0$,
one has
\begin{equation}\label{N-K-B}\log N(K,tB)\ls \frac{c_1nI_1(K,B)}{t}+\log 2,\end{equation}
where $c_1>0$ is an absolute constant.
\end{lemma}

\noindent {\it Proof.} We define a Borel probability measure
on ${\mathbb R}^n$ by
\begin{equation}\mu(A)=\frac{1}{c_K}\int_Ae^{-p_K(x)}dx,\end{equation}
where $p_K$ is the Minkowski functional of $K$ and
$c_K=\int_{{\mathbb R}^n}\exp (-p_K(x))dx$. A simple computation,
based on the fact that $\{x\in {\mathbb R}^n:p_K(x)\ls t\} = tK$
for any $t>0$, shows that $c_K=n!$.

Let $\{x_1,\ldots ,x_N\}$ be a subset of $K$ which is maximal with
respect to the condition $\|x_i-x_j\|_B\gr t$ for $i\neq j$. Then
$K\subseteq\bigcup_{i\ls N}(x_i+tB)$, and hence $N(K,tB)\ls N$. Let
$a>0$. Note that if we set $y_i=(2a/t)x_i$, by the subadditivity of
$p_K$ and the fact that $p_K(x_i)\ls 1$, we have
\begin{equation} \mu (y_i+aB) \gr\frac{1}{c_K}\int_{aB}e^{-p_K(x)}e^{-p_K(y_i)}dx\gr e^{-2a/t}\mu (aB).\end{equation}
The bodies $y_i+aB$ have disjoint interiors, therefore $Ne^{-2a/t}\mu (aB)\ls 1$.
It follows that
\begin{equation}N(K,tB)\ls 2e^{2a/t}(\mu (aB))^{-1}.\end{equation}
Now, we choose $a>0$ so that $\mu (aB)\gr 1/2$. A simple computation
shows that
\begin{equation}J := \int_{{\mathbb R}^n}\|x\|_K\, d\mu(x)=(n+1)I_1(K,B).\end{equation}
By Markov's inequality, $\mu(2JB)\gr 1/2$, so if we choose $a=2J$, we get
\begin{equation}N(K,tB) \ls  2\exp (4J/t)\ls 2\exp\bigl(4(n+1)I_1(K,B)/t\bigr)\end{equation}
for every $t>0$. $\hfill\Box$

\smallskip

\begin{remark}\label{remark:isotropic-covering}
\rm (i) In the case that $B$ is the Euclidean ball $B_2^n$ and $K$
is an isotropic convex body, we have that $I_1(K,B)\ls \sqrt{n}L_K$
and therefore
\begin{equation}\label{isotropic-covering1} \log N(K,tB_2^n)\ls \frac{c'_1n^{3/2}L_K}{t}\end{equation}
for any $t>0$ (for very large $t$ the estimate is trivially true, since every isotropic
body $K$ satisfies the inclusion $K\subseteq cnL_K B_2^n$ for some absolute constant $c$).
Given (\ref{isotropic-condition2}), this is essentially the best way we
can apply Lemma \ref{covering-K-by-tB} when $B=B_2^n$. This version of the lemma
appeared in the Ph.D. Thesis of Hartzoulaki \cite{Har}. The idea of using
$I_1(K, B_2^n)$ as a parameter in entropy estimates for isotropic convex bodies
comes from \cite{MP1}. It was also used in \cite{LMP} for a proof of the low
$M^{\ast}$-estimate in the case of quasi-convex bodies.

\smallskip

\noindent (ii) Knowing that we have for any set $S$,
\begin{equation}N(S-S, 2B_2^n)= N(S-S,B_2^n-B_2^n)\ls N(S,B_2^n)^2,\end{equation}
we can use (\ref{isotropic-covering1}) to also get an upper bound for the covering
numbers of the difference body of an isotropic convex body $K$ by the Euclidean ball:
\begin{equation}\log N(K-K,tB_2^n)\ls \frac{2c'_1n^{3/2}L_K}{t}.\end{equation}

\noindent (iii) Lemma \ref{covering-K-by-tB} is also related to the problem
of estimating the mean width of an isotropic convex body $K$, namely the parameter
$w(K):=\int_{S^{n-1}}h_K(\theta)d\sigma(\theta)$ where $h_K$ is the support function of $K$
and $\sigma$ is the uniform probability measure on $S^{n-1}$. The best upper
bound we have is $w(K)\ls cn^{3/4}L_K$ (there are several arguments leading to this estimate;
see \cite{GPV} and the references therein). It is known (see e.g. \cite[Theorem 5.6]{GM})
that an improvement of the form
\begin{equation}\log N(K,tB_2^n)\ls \frac{c'_1n^{3/2}L_K}{t^{1+\delta }}\end{equation}
(for some $\delta >0$) in (\ref{isotropic-covering1}) would immediately imply a better
bound for $w(K)$ in the isotropic case.
\end{remark}

\smallskip

The next lemma allows us to bound the dual covering numbers $N(B_2^n, tK^{\circ})$.

\begin{lemma}\label{isotropic-covering2}
Let $K$ be a convex body in ${\mathbb R}^n$ which contains $0$ in
its interior. For every $t>0$ we set $A(t):=t\log N(K,tB_2^n)$ and
$B(t):=t\log N(B_2^n,tK^\circ)$. Then, one has
\begin{equation}\sup_{t>0}B(t)\ls 16\sup_{t>0}A(t).\end{equation}
In particular, if $K$ is isotropic $($or a translate of an
isotropic convex body which still contains $0$ in its interior$)$,
then
\begin{equation}\label{N-B-K}\log N(B_2^n,tK^\circ)
\ls\log N\bigl(B_2^n,t(K-K)^\circ\bigr)\ls
\frac{c_2n^{3/2}L_K}{t},\end{equation} where $c_2>0$ is an absolute
constant.
\end{lemma}

\noindent {\it Proof.} We use a well-known idea from \cite{TJ} (see
also \cite[Section 3.3]{LT}). For any $t>0$ we have
$(t^2K^\circ)\cap (4K)\subseteq 2tB_2^n$. Passing to the polar
bodies we see that
\begin{equation}B_2^n\subseteq {\rm conv}\left ( \frac{t}{2}K^{\circ
},\frac{2}{t}K\right )\subseteq \frac{t}{2}K^{\circ
}+\frac{2}{t}K.\end{equation} We write
\begin{align}
N(B_2^n,tK^\circ)&\ls N\left (\frac{t}{2}K^{\circ }+\frac{2}{t}K,tK^{\circ }\right )
= N\left (\frac{2}{t}K,\frac{t}{2}K^{\circ }\right )\\
\nonumber &\ls N\left (\frac{2}{t}K,\frac{1}{4}B_2^n\right )N\left
(\frac{1}{4}B_2^n,\frac{t}{2}K^{\circ }\right )\\
\nonumber &= N\left (K,\frac{t}{8}B_2^n\right )N(B_2^n,2tK^{\circ
}).
\end{align} Taking logarithms we get \begin{equation}
B(t) \ls 8 A(t/8)+\frac{1}{2}B(2t),
\end{equation} for all $t>0$. This implies that \begin{equation}
B:=\sup_{t>0}B(t)\ls 16 A,
\end{equation} and the result follows. $\hfill\Box $

\medskip

The last covering lemma is from \cite{iso-sym} and shows that the
volume $|{\rm conv}\bigl (K\cup L\bigr )|$ of the convex hull of two
convex bodies $K$ and $L$ is essentially bounded by $N(L,K)\,|K|$,
provided that $L\subseteq bK$ for some ``reasonable" $b\gr 1$.

\begin{lemma}\label{convex-hull}
Let $L$ be a convex body and let $K$ be a symmetric convex body in
${\mathbb R}^n$. Assume that $L\subseteq bK$ for some $b\gr 1$. Then
\begin{equation}\big|{\rm conv}\bigl(K \cup L\bigr)\big| \ls 3enb\ N(L,K)|K|.\end{equation}
\end{lemma}

\noindent {\it Proof.} By the definition of $N\equiv N(L,K)$, there
exist $x_{1},\ldots,x_{N}\in {\mathbb R}^n$ such that $(x_i+K)\cap
L\neq\emptyset$ for every $i=1,\ldots,N$, and
\begin{equation}L\subseteq \bigcup_{i=1}^{N}(x_{i}+K).\end{equation}
From the symmetry of $K$ and the fact that $L \subseteq bK$, it
follows that, for every $i=1,\ldots,N$,
\begin{equation}\label{eqp:CovLem1}x_i\in L+K\subseteq (1+b)K.\end{equation}
Now, for every $\alpha, \beta\in [0,1]$ with $\alpha +\beta =1$, we
have that
\begin{align}
\alpha L + \beta K &\subseteq
\bigcup_{i=1}^{N}(\alpha x_{i}+\alpha K)+\beta K
= \bigcup_{i=1}^{N}\bigl(\alpha x_{i}+(\alpha+\beta)K\bigr)\\
\nonumber &= \bigcup_{i=1}^{N}(\alpha x_{i}+K),
\end{align}
and therefore
\begin{equation}\label{eqp:CovLem2}
{\rm conv}(L\cup K)=\! \bigcup_{0\ls\alpha\ls 1}\bigl(\alpha L + (1-\alpha)K\bigr)\subseteq
\bigcup_{i=1}^{N}\bigcup_{0\ls\alpha\ls 1}(\alpha x_{i}+K).
\end{equation}
We set $T=2n$ and consider $\lceil bT \rceil$ numbers $\alpha_{j}$
equidistributed in $[0,1]$, $j=1,\ldots,\lceil bT \rceil$. From
(\ref{eqp:CovLem1}) and (\ref{eqp:CovLem2}) it follows that: for
every $z\in {\rm conv}(L\cup K)$ there exist $\alpha, \alpha_{j}\in
[0,1]$, with distance $|\alpha -\alpha_{j}|\ls \frac{1}{bT}$, such
that
\begin{equation}\label{eqp:CovLem3}
z\in \alpha x_{i}+K= \alpha_{j}x_{i}+(\alpha -\alpha_{j})x_{i}+K\subseteq \alpha_{j}x_{i}+\left(\frac{1+ b}{bT}+1\right)K.
\end{equation}
We observe that
\begin{equation}\frac{1+ b}{bT}=\frac{1+b}{2nb}\ls\frac{1}{n}\end{equation}
because $b\gr 1$ and $\frac{1+b}{2}\ls b$, so (\ref{eqp:CovLem3})
gives us that
\begin{equation}z\in \alpha_{j}x_{i}+\left(1+\frac{1}{n}\right)K.\end{equation}
Going back to (\ref{eqp:CovLem2}), we see that
\begin{equation}{\rm conv}(L\cup K)\subseteq \bigcup_{i=1}^{N}\bigcup_{j=1}^{\lceil bT \rceil}
\left\{\alpha_{j}x_{i}+\left(1+\frac{1}{n}\right)K\right\}.\end{equation}
Then,
\begin{align}
|{\rm conv}(L\cup K)| &\ls N\lceil bT \rceil\left(1+\frac{1}{n}\right)^{n}|K|\ls \tfrac{3}{2}bTeN|K|\\
\nonumber &= 3enb\,N\!(L,K)|K|,
\end{align}
which is our claim. $\hfill\Box$

\medskip

\noindent {\bf 2.3. The method of convex perturbations}. In \cite{perturbation}
Klartag gave an affirmative answer to the following question: even if we don't
know that every convex body in ${\mathbb R}^n$ has bounded isotropic constant,
given a body $K$ can we find a second body $T$ ``geometrically close'' to $K$
with isotropic constant $L_T\simeq 1$? Here when we say that $K$ and $T$ are
``geometrically close'', we will mean that there exists an absolute constant
$c>0$ such that for some $x,y\in {\mathbb R}^n$,
\begin{equation}\label{GeometricDistance1}\frac{1}{c}(T-x)\subseteq K-y\subseteq c(T-x).\end{equation}
The method Klartag used is based on two key observations. The first one is that
in order to find a body $T$ close to $K$ which has bounded isotropic constant,
it suffices to define a positive log-concave function on $K$ (vanishing everywhere else)
with bounded isotropic constant and the extra property that its range is not too large.

\begin{proposition}\label{FunctionsToBodies}
Let $K$ be a convex body in ${\mathbb R}^n$ and let $f:K\to
(0,\infty )$ be a log-concave function such that
\begin{equation}\sup_{x\in K}f(x)\ls m^n\inf_{x\in
K}f(x)\end{equation} for some $m>1$. Let $x_0$ be the barycenter of $f$,
i.e. $x_0=\int_{{\mathbb R}^n} xf(x)\,dx/\int_{{\mathbb R}^n}f(x)\,dx$,
and set $g(x)= f(x+x_0)$. Then, for the centered convex body $T:=K_{n+1}(g)$,
defined as in $(\ref{Kp-definition})$, we have that $L_f\simeq L_T$ and
\begin{equation}\frac{1}{m}T\subseteq K-x_0\subseteq mT.\end{equation}
\end{proposition}

The second observation is that a family of suitable candidates for the
function $f$ we need so as to apply Proposition \ref{FunctionsToBodies} can
be found through the logarithmic Laplace transform on $K$. In general, the
logarithmic Laplace transform of a finite Borel measure $\mu $ on
${\mathbb R}^n$ is defined by
\begin{equation}\Lambda_{\mu }(\xi )=\log\left (\int_{{\mathbb R}^n}e^{\langle\xi
,x\rangle }\frac{d\mu (x)}{\mu({\mathbb R}^n)}\right ).\end{equation}
In \cite{perturbation}, Klartag makes use of the following properties of $\Lambda_{\mu }$:

\begin{proposition}\label{LaplaceTransformProperties}
Let $\mu =\mu_K$ denote the Lebesgue measure on some convex body $K$ in ${\mathbb R}^n$. Then,
\begin{equation}\bigl(\nabla\Lambda_{\mu }\bigr)({\mathbb R}^n)={\rm
int}(K)\end{equation} $($actually, for the arguments in {\rm
\cite{perturbation}} and for our proof here, it suffices to know
that $\bigl(\nabla\Lambda_{\mu }\bigr)({\mathbb R}^n)\subseteq K)$.
If $\mu_{\xi }$ is the probability measure on ${\mathbb R}^n$ with
density proportional to the function $e^{\langle \xi ,x\rangle }{\bf
1}_K(x)$, then
\begin{equation}b(\mu_{\xi })=\nabla \Lambda_{\mu }(\xi )\quad {\rm and}\quad {\rm
Hess}\,(\Lambda_{\mu }(\xi ))={\rm Cov}(\mu_{\xi }).\end{equation}
Moreover, the map $\nabla \Lambda_{\mu }$, which is one-to-one, transports the measure
$\nu $ with density $\det {\rm Hess}\,(\Lambda_{\mu })$ to
$\mu $. In other words, for every continuous non-negative function
$\phi :{\mathbb R}^n\to {\mathbb R}$,
\begin{equation}\int_K\phi (x)\,dx=\int_{{\mathbb R}^n}\phi (\nabla\Lambda_{\mu
}(\xi ))\det {\rm Hess}(\Lambda_{\mu }(\xi ))\,d\xi =\int_{{\mathbb R}^n}\phi
(\nabla\Lambda_{\mu }(\xi ))d\nu (\xi ).\end{equation}
\end{proposition}

\smallskip

Klartag's approach has been recently applied in \cite{DPV} where
Dadush, Peikert and Vempala provide an algorithm for enumerating
lattice points in a convex body, with applications to integer
programming and problems about lattice points. They use the
techniques of \cite{perturbation} in order to give an expected
$2^{O(n)}$-time algorithm for computing an $M$-ellipsoid for any
convex body in ${\mathbb R}^n$.

\section{Proof of the reverse Santal\'{o} inequality}

We now prove the reverse Santal\'{o} inequality using the results
that were described in Section 2. The proof consists of three steps
which roughly are the following: (i) we obtain a lower bound for
the volume product $s(K)$ which is optimal up to the value of the
isotropic constant $L_K$ of $K$, (ii) by adapting Klartag's main
argument from \cite{perturbation} we show that every symmetric convex
body $K$ has bounded geometric distance (in the sense defined in
(\ref{GeometricDistance1})) from a second convex body $T$ whose isotropic
constant $L_T$ can be expressed in terms of $s(K)$, and (iii) we use the
lower bound for $s(T)$ in terms of $L_T$, and the fact that $s(K)$ and $s(T)$
are comparable, to get a lower bound for $s(K)$ in which $L_K$ does not appear anymore.

\medskip

\noindent {\bf 3.1. Lower bound involving the isotropic constant}.
Our first step will be to prove the following lower bound for
$s(K)$.

\begin{proposition}\label{reverse-santalo-with-LK}
Let $K$ be a convex body in $\mathbb R^n$ which contains $0$ in its
interior. Then
\begin{equation}4|K|^{1/n}|nK^{\circ}|^{1/n}\gr |K-K|^{1/n}
|n(K-K)^{\circ}|^{1/n}\gr\frac{c_1}{L_K},
\end{equation}where $c_1>0$ is an absolute constant.
\end{proposition}

\noindent {\it Proof.} We may assume that $|K|=1$. From the
Brunn-Minkowski inequality and the classical Rogers--Shephard
inequality (see \cite{RS}), we have $2 \ls|K-K|^{1/n}\ls 4$. Since
$(K-K)^{\circ}\subseteq K^{\circ}$, we immediately see that
\begin{equation}|K|^{1/n}|nK^{\circ}|^{1/n}\gr
\frac{1}{4}|K-K|^{1/n}|n(K-K)^{\circ}|^{1/n},\end{equation}
so it remains to prove the second inequality. Since
\begin{equation}\big|T(K)-T(K)\big|\big|\bigl(T(K)-T(K)\bigr)^{\circ}\big|
=|K-K||(K-K)^{\circ}|\end{equation} for any invertible affine
transformation $T$ of $K$, we may assume for the rest of the proof that $K$ is
isotropic. We define
\begin{equation}K_1:=\frac{K-K}{L_K}\cap\overline{B}_2^n\end{equation}
and observe that the inclusion $K_1\subseteq \overline{B}_2^n$ implies that
$\overline{B}_2^n\subseteq c_1nK_1^{\circ}$ for some absolute constant $c_1$.
Moreover,
\begin{equation}nK_1^{\circ}\simeq {\rm conv}\{ nL_K(K-K)^{\circ},\overline{B}_2^n\},\end{equation}
therefore we can apply Lemma \ref{convex-hull} with $L=\overline{B}_2^n$ and
$K=nL_K(K-K)^{\circ }$ to bound $|nK_1^{\circ}|$ from above; note that in this
case $b\simeq\sqrt{n}$, because $K-K\subseteq cnL_KB_2^n$ since we have assumed $K$
isotropic (see \cite[Theorem 1.2.4]{Gian}), and hence
$\overline{B}_2^n\subseteq c'\sqrt{n}\bigl(nL_K(K-K)^{\circ}\bigr)$ for some
absolute constants $c, c'$. Using also (\ref{N-B-K})
from Lemma \ref{isotropic-covering2} (with $t\simeq\sqrt{n}L_K$), we see that
\begin{align}
c_1^{-n} &\ls |n K_1^{\circ}| \ls c_2^n|{\rm conv} \{ nL_K(K-K)^{\circ}, \overline{B}_2^n\}| \\
\nonumber &\ls c_3^nn^{3/2} |nL_K(K-K)^{\circ}|\,N\!\left(\overline{B}_2^n, nL_K(K-K)^{\circ}\right) \\
\nonumber &\ls c_3^nn^{3/2} |nL_K (K-K)^{\circ}|\,N\!\left(B_2^n, c_4\sqrt{n}L_K(K-K)^{\circ}\right)\\
\nonumber &\ls e^{c_5 n}|nL_K (K-K)^{\circ}|.
\end{align}
This shows that there exists an absolute constant $c'_1$ so that
\begin{equation}|nL_K(K-K)^{\circ}|^{1/n}\gr c'_1,\end{equation}
and since $|K-K|^{1/n}\gr 2$, we have proven that
\begin{equation}\hspace{2,7cm} |K-K|^{1/n}|(K-K)^{\circ}|^{1/n}\gr \frac{2c'_1}{nL_K}. \hspace{3,3cm}\Box\end{equation}

\medskip

\noindent {\bf 3.2. A variant of Klartag's argument}. Our second
step will be to show that every convex body $K$ in ${\mathbb R}^n$ has
bounded geometric distance from a second convex body $T$ whose isotropic
constant $L_T$ can be bounded in terms of $s(K\!-\!K)$.

\begin{proposition}\label{klartag-argument}
Let $K$ be a convex body in ${\mathbb R}^n$. For every
$\varepsilon\in (0,1)$ there exist a centered convex body
$T\subset {\mathbb R}^n$ and a point $x\in {\mathbb R}^n$ such that
\begin{equation}\label{klartag1}
\frac{1}{1+\varepsilon}T\subseteq K+x \subseteq (1+\varepsilon)T
\end{equation}
and
\begin{equation}\label{klartag2}
L_T\ls\frac{c_2}{\sqrt{\varepsilon ns(K\!-\!K)^{1/n}}},
\end{equation}
where $c_2>0$ is an absolute constant.
\end{proposition}

\noindent{\it Proof.} We may assume that $K$ is centered and that $|K\!-\!K|=1$. Indeed, once we
prove the proposition for $\tilde{K}:= (K\!-\!{\rm bar}(K))/|K\!-\!K|^{1/n}$ and some
$\varepsilon\in (0,1)$, and find a convex body $T$ which satisfies (\ref{klartag1}) and
(\ref{klartag2}) with $\tilde{K}$ instead of $K$, it will immediately hold that the pair
$(K,|K\!-\!K|^{1/n}T)$ also satisfies these properties, because $L_T$ and $s(K\!-\!K)$ are
affine invariants.

Recall from Proposition \ref{LaplaceTransformProperties} that if $\mu=\mu_K$ is
the Lebesgue measure restricted on $K$, then the function
$\nabla \Lambda_{\mu }$ transports the measure $\nu $ with density
\begin{equation}\frac{d\nu}{d\xi}=\det {\rm Hess}\,(\Lambda_{\mu }(\xi))\equiv\det {\rm Cov}(\mu_{\xi })\end{equation}
to $\mu$. This implies that
\begin{equation}\nu({\mathbb R}^n)=\int_{{\mathbb R}^n}{\bf 1}\det {\rm Hess}\,(\Lambda_{\mu }(\xi))\,d\xi=
\int_K{\bf 1}\,dx = |K|\ls|K\!-\!K|=1.\end{equation}
Thus, for every $\varepsilon >0$ we may write
\begin{multline}
|\varepsilon n(K\!-\!K)^{\circ }|\min_{\xi\in\varepsilon n(K\!-\!K)^{\circ }}\det
{\rm Cov}(\mu_{\xi })\ls
\\
\ls  \int_{\varepsilon n(K\!-\!K)^{\circ }}\det {\rm Cov}(\mu_{\xi })\,d\xi  =\nu (\varepsilon n(K\!-\!K)^{\circ })\ls 1,
\end{multline}
which means that there exists $\xi\in\varepsilon n(K\!-\!K)^{\circ }$ such
that
\begin{equation}\det {\rm Cov}(\mu_{\xi })=\min_{\xi'\in\varepsilon n(K\!-\!K)^{\circ }}\det {\rm Cov}(\mu_{\xi' })
\ls |\varepsilon n (K\!-\!K)^{\circ}|^{-1}=\bigl(\varepsilon ns(K\!-\!K)^{1/n}\bigr)^{-n}\end{equation}
(where the last equality holds because $|K-K|=1$). Now, from the definition of $\mu_{\xi}$ and (\ref{definition-isotropic})
we have that
\begin{equation}L_{\mu_{\xi }}=\left (\frac{\sup_{x\in K}e^{\langle\xi ,x\rangle }}
{\int_Ke^{\langle\xi ,x\rangle }dx}\right )^{_{\frac{1}{n}}}[\det {\rm
Cov}(\mu_{\xi } )]^{\frac{1}{2n}}.\end{equation} Since
$\xi\in\varepsilon n(K\!-\!K)^{\circ }$ and $K\cup(-K)\subset K\!-\!K$, we know that
$|\langle\xi ,x\rangle| \ls\varepsilon n$ for all $x\in K$, therefore
$\sup_{x\in K}e^{\langle\xi ,x\rangle} \ls \exp(\varepsilon n)$.
On the other hand, since $K$ is centered, from Jensen's inequality we have that
\begin{equation}\frac{1}{|K|}\int_K e^{\langle\xi ,x\rangle }dx\gr\exp \left (\frac{1}{|K|}\int_K\langle\xi
,x\rangle \,dx\right )=1,\end{equation}
which means that $\int_Ke^{\langle\xi ,x\rangle }dx\gr |K|\gr 4^{-n}|K\!-\!K|$ by the
Rogers-Shephard inequality. Combining all these we get
\begin{equation}L_{\mu_{\xi }}\ls \frac{4e^{\varepsilon
}}{\sqrt{\varepsilon ns(K\!-\!K)^{1/n}}}.\end{equation}
\smallskip

\noindent Finally, we note that the function $f_{\xi
}(x)=e^{\langle \xi ,x\rangle }{\bf 1}_K(x)$ (which is proportional
to the density of $\mu_{\xi}$) is obviously log-concave and satisfies
\begin{equation}
\sup_{x\in {\rm supp}(f_{\xi})}f_{\xi }(x)\ls e^{2\varepsilon
n}\inf_{x\in {\rm supp}(f_{\xi})}f_{\xi }(x)
\end{equation}
(since $|\langle\xi ,x\rangle| \ls\varepsilon n$ for all $x\in K$).
Therefore, applying Proposition \ref{FunctionsToBodies}, we can
find a centered convex body $T_{\xi}$ in ${\mathbb R}^n$ such that
\begin{equation}L_{T_{\xi}}\simeq L_{f_{\xi}}=L_{\mu_{\xi }}\ls
\frac{4e^{\varepsilon}}{\sqrt{\varepsilon ns(K\!-\!K)^{1/n}}}\end{equation}
and
\begin{equation}\frac{1}{e^{2\varepsilon }}T_{\xi}\subseteq K-b_{\xi}\subseteq e^{2\varepsilon}T_{\xi}\end{equation}
where $b_{\xi}$ is the barycenter of $f_{\xi}$. Since $e^{2\varepsilon }\ls
1+c\:\!\varepsilon $ when $\varepsilon\in (0,1)$, the result follows. $\hfill\Box $

\medskip

\noindent {\bf 3.3. Removing the isotropic constant}. Combining the
previous two results we can remove the isotropic constant $L_K$ from
the lower bound for $s(K)^{1/n}$.

\begin{theorem}\label{reverse-santalo}Let $K$ be a convex body in $\mathbb R^n$ which contains $0$ in its
interior. Then
\begin{equation} |K|^{1/n}|nK^{\circ}|^{1/n}\gr c_3,\end{equation}
where $c_3>0$ is an absolute constant.
\end{theorem}

\noindent {\it Proof.} Since $|K|^{1/n}|nK^{\circ}|^{1/n}\gr
\frac{1}{4}|K-K|^{1/n} |n(K-K)^{\circ}|^{1/n}$, we may assume
for the rest of the proof that $K$ is symmetric. Using Proposition
\ref{klartag-argument} with $\varepsilon =1/2$, we find a convex
body $T\subset {\mathbb R}^n$ and a point $x\in {\mathbb R}^n$ such that
\begin{equation}\label{klartag3}
\frac{2}{3}T\subseteq K+x \subseteq \frac{3}{2}T
\end{equation}
and $L_T\ls c_0/\sqrt{ns(K)^{1/n}}$ for some absolute constant $c_0>0$.
Proposition \ref{reverse-santalo-with-LK} shows that
\begin{equation}|T-T|^{1/n} |n(T-T)^{\circ}|^{1/n}\gr
\frac{c_1}{L_T},\end{equation} where $c_1>0$ is an absolute
constant too. Observe that $\frac{2}{3}(T-T)\subseteq
K-K=2K\subseteq\frac{3}{2}(T-T)$, and thus
$K^{\circ}\supseteq\frac{4}{3}(T-T)^{\circ }$. Therefore, combining the above, we get
\begin{align}ns(K)^{1/n} &=|nK^{\circ}|^{1/n}|K|^{1/n} \gr \frac{4}{9} |n(T-T)^{\circ}|^{1/n}|T-T|^{1/n}
\\ \nonumber
&\gr \frac{c'_1}{L_T}\gr c_2\sqrt{ns(K)^{1/n}},\end{align} and so it follows that
\begin{equation}s(K)^{1/n}\gr\frac{c_3}{n}\end{equation}
with $c_3=c_2^2$. This completes the proof. $\hfill\Box $

\medskip

Having proved the reverse Santal\'{o} inequality, one can go back to
Proposition \ref{klartag-argument} and insert the lower bound for
$s(K\!-\!K)$. This is the last step in Klartag's solution of the
isomorphic slicing problem.

\begin{theorem}[Klartag]\label{klartag-isomorphic}
Let $K$ be a convex body in ${\mathbb R}^n$. For every
$\varepsilon\in (0,1)$ there exist a centered convex body
$T\subset{\mathbb R}^n$ and a point $x\in {\mathbb R}^n$ such that
\begin{equation}\label{klartag5}
\frac{1}{1+\varepsilon}T\subseteq K+x \subseteq (1+\varepsilon)T
\end{equation}
and
\begin{equation}\label{klartag6}
L_T\ls\frac{c_4}{\sqrt{\varepsilon }},
\end{equation}
where $c_4>0$ is an absolute constant.
\end{theorem}

\section{$M$-ellipsoids and the reverse Brunn-Minkowski inequality}

We can now prove the existence of $M$-ellipsoids for any convex body
and, as a consequence, the reverse Brunn--Minkowski inequality.

\medskip

\noindent {\bf 4.1. Existence of $M$-ellipsoids}. Let $K$ be a
centered convex body in ${\mathbb R}^n$. We will give a proof of
the existence of an $M$-ellipsoid for $K$. The next Proposition is
the first step.

\begin{proposition}\label{existence-of-EK}
Let $K$ be a centered convex body in $\mathbb R^n$. Then there
exists an ellipsoid ${\cal{E}}_K$ such that $|K|=|{\cal{E}}_K|$ and
\begin{equation}\max\{\log{ N(K,t{\cal{E}}_K)},\log{ N({\cal{E}}_K^{\circ},tK^{\circ})}\}
\ls \frac{cn}{t}\end{equation}for all $t>0$, where $c>0$ is an
absolute constant.
\end{proposition}

\noindent {\it Proof.} Applying Proposition \ref{klartag-isomorphic}, we can find
a centered convex body $T$ with isotropic constant $L_T\ls C$ such that
\begin{equation}\frac{2}{3}T\subseteq K+x\subseteq \frac{3}{2}T\end{equation}
for some $x\in {\mathbb R}^n$. Let $Q(T)$ be an isotropic position of $T$. From Remark
\ref{remark:isotropic-covering}(ii) and Lemma \ref{isotropic-covering2} we know that
\begin{equation}\label{ball-covering1}\max\{ \log{ N\bigl(Q(T)-Q(T),t\sqrt{n}B_2^n\bigr)},
\log{ N\bigl(B_2^n,t\sqrt{n}(Q(T)-Q(T))^{\circ}\bigr)}\}\ls \frac{cn}{t}\end{equation}
for every $t>0$. Since
\begin{equation}\frac{2}{3}(Q(T)-Q(T))\subseteq Q(K)-Q(K)\subseteq \frac{3}{2}(Q(T)-Q(T))\end{equation}
and $Q(K)\subseteq Q(K)-Q(K),\; (Q(K)-Q(K))^{\circ}\subseteq (Q(K))^{\circ}$, from (\ref{ball-covering1}) it follows that
\begin{equation}\label{ball-covering2}\max\{ \log{ N\bigl(Q(K),t\sqrt{n}B_2^n\bigr)},
\log{ N\bigl(B_2^n,t\sqrt{n}(Q(K))^{\circ}\bigr)}\}\ls
\frac{c'n}{t}\end{equation} for every $t>0$. We define
${\cal{E}}_K:= Q^{-1}(a\sqrt{n}B_2^n)$ where $a$ is chosen so that
$|Q(K)|=|a\sqrt{n}B_2^n|$ (equivalently, so that
$|{\cal{E}}_K|=|K|$), and from (\ref{ball-covering1}) we get that
\begin{equation}\max\{\log{ N(K,t{\cal{E}}_K)},\log{ N({\cal{E}}_K^{\circ},tK^{\circ})}\}
\ls \frac{c'an}{t}\end{equation}
for all $t>0$. It remains to observe that
\begin{equation}|\sqrt{n}B_2^n|^{1/n}\simeq 1=|Q(T)|^{1/n}\simeq |Q(K+x)|^{1/n} = |Q(K)|^{1/n},\end{equation}
whence it follows that $a\simeq 1$. $\hfill\Box $

\medskip

We now recall some standard entropy estimates which are valid for
arbitrary convex bodies in ${\mathbb R}^n$.

\begin{lemma}\label{lemma-KM2}
Let $K$ and $L$ be convex bodies in ${\mathbb R}^n$. If $L$ is
symmetric, then
\begin{equation}\label{eq:entropy-estimates1}N(K,L) \ls \frac{|K + L/2|}{|L/2|}\ls 2^n\frac{|K+L|}{|L|},\end{equation}
whereas in the general case
\begin{equation}\label{eq:entropy-estimates2}N(K,L)\ls 4^n\frac{|K+L|}{|L|}.\end{equation}
Moreover,
\begin{equation}\label{eq:entropy-estimates3}\frac{|K + L|}{|L|} \ls 2^n N(K,L).\end{equation}
\end{lemma}

\noindent {\it Proof.} The proof of (\ref{eq:entropy-estimates3}) is
an easy consequence of the definitions. To prove (\ref{eq:entropy-estimates1}),
note that if $N$ is a maximal subset of $K$ with respect to the property
\begin{equation}x,y \in N \ \hbox{and}\ x\neq y \Rightarrow \|x-y\|_L\gr 1,\end{equation}
then $K\subseteq \bigcup_{x\in N}(x+L)$, while every two sets $x+L/2, y+L/2 \ (x,y\in N)$
have disjoint interiors when $x\neq y$.

Finally, when $L$ is not necessarily symmetric, we recall that $N(K+x,L+y)=N(K,L)$
for every $x,y\in {\mathbb R}^n$, and also that the ratio $|K+L|/|L|$ obviously
remains unaltered if we translate $K$ or $L$. Hence, we can assume that $L$ is centered,
in which case it follows from \cite[Corollary 3]{MP2} that
\begin{equation}|L\cap (-L)|\gr 2^{-n}|L|.\end{equation}
But then, from (\ref{eq:entropy-estimates1}) we get that
\begin{equation}N(K,L)\ls N(K, L\cap (-L))\ls 2^n\frac{|K + (L\cap (-L))|}{|L\cap (-L)|}\ls
4^n\frac{|K + L|}{|L|},\end{equation} and we have
(\ref{eq:entropy-estimates2}). $\hfill\Box $

\begin{corollary}\label{corollary-KM2}
Let $K$ and $L$ be two convex bodies in ${\mathbb R}^n$. Then,
\begin{equation}N(K,L)^{1/n}\simeq \frac{|K + L|^{1/n}}{|L|^{1/n}}.\end{equation}
It also follows that if $K$ and $L$ have the same volume, then
\begin{equation}N(K,L)^{1/n}\ls 8N(L,K)^{1/n}.\end{equation}
\end{corollary}

\medskip

Combining Proposition \ref{existence-of-EK} with the classical
Santal\'{o} inequality and Corollary \ref{corollary-KM2},
we can now prove the existence of $M$-ellipsoids
for any centered convex body in ${\mathbb R}^n$.

\begin{theorem}\label{M-ellipsoids}
Let $K$ be a centered convex body in $\mathbb R^n$. There exists an
ellipsoid ${\cal E}_K$ such that $|K|=|{\cal E}_K|$ and
\begin{equation}\max\big\{\log N(K,{\cal E}_K), \log N({\cal E}_K,K), \log N( K^{\circ},{\cal E}_K^{\circ}),
\log N({\cal E}_K^{\circ},K^{\circ})\big\}\ls cn,\end{equation}
where $c>0$ is an absolute constant.
\end{theorem}

\noindent {\it Proof.} Let ${\cal E}_K$ be the ellipsoid defined in
Proposition \ref{existence-of-EK}. It immediately follows that
\begin{equation}\max\big\{N(K,{\cal E}_K), N({\cal E}_K^{\circ},K^{\circ})\big\}\ls \exp(cn).\end{equation}
For the other two covering numbers we use Lemma \ref{lemma-KM2}:
$N({\cal E}_K,K)\ls 8^n N(K,{\cal E}_K)$, which means that $\log N({\cal E}_K,K)\ls
(\log 8)n + \log N(K,{\cal E}_K)$. Similarly,
\begin{equation}N( K^{\circ},{\cal E}_K^{\circ})\ls
2^n\frac{|K^{\circ} + {\cal E}_K^{\circ}|}{|{\cal E}_K^{\circ}|}\ls
2^n\frac{|K^{\circ} + {\cal E}_K^{\circ}|}{|K^{\circ}|}\ls 4^n N({\cal E}_K^{\circ},K^{\circ}),
\end{equation}
where we have also used the fact that $|K|=|{\cal E}_K|\Rightarrow |K^{\circ}|\ls |{\cal E}_K^{\circ}|$
from the classical Santal\'{o} inequality. This completes the proof.
$\hfill\Box $

\medskip

\noindent {\bf 4.2. Reverse Brunn--Minkowski inequality}. As a
consequence of Theorem \ref{M-ellipsoids} and Corollary
\ref{corollary-KM2}, we get the ``reverse" Brunn-Minkowski
inequality.

\begin{theorem}\label{theorem-ReverseBM}
Let $K$ be a centered convex body in $\mathbb R^n$. There exists an
ellipsoid ${\cal E}_K$ such that $|K|=|{\cal E}_K|$ and for every convex body
$T$ in ${\mathbb R}^n$,
\begin{align}
\label{eq:towardsReverseBM1} e^{-(c+\log\!8)}\;|{\cal
E}_K+T|^{1/n}&\ls |K + T|^{1/n}\ls e^{c+\log\!8}\;|{\cal E}_K +
T|^{1/n},
\\
\label{eq:towardsReverseBM2} e^{-(c+\log\!8)}\;|{\cal E}_K^{\circ
}+T|^{1/n}&\ls |K^{\circ}+T|^{1/n}\ls e^{c+\log\!8}\;|{\cal
E}_K^{\circ}+T|^{1/n},
\end{align}
where $c$ is the constant we found in Theorem {\rm \ref{M-ellipsoids}}.
\end{theorem}

\noindent {\it Proof.} Let ${\cal E}_K$ be the ellipsoid defined in
Proposition \ref{existence-of-EK}. Using Lemma \ref{lemma-KM2}, we
can write
\begin{align}
|{\cal E}_K+T|^{1/n}&\ls 2|T|^{1/n}N({\cal E}_K,T)^{1/n}\ls 2|T|^{1/n}N({\cal E}_K,K)^{1/n}N(K,T)^{1/n}
\\
\nonumber
&\ls 2e^c |T|^{1/n}N(K,T)^{1/n}\ls 8e^c |K+T|^{1/n}.
\end{align}
The same reasoning gives us the second part of
(\ref{eq:towardsReverseBM1}) and (\ref{eq:towardsReverseBM2}).
$\hfill\Box $

\medskip

\begin{remark}
\rm We usually say that a centered convex body $K$ is in
$M$-position if the ellipsoid ${\cal E}_K$ that we look for in
Theorem \ref{M-ellipsoids} can be taken to be a multiple of the
Euclidean ball. Obviously, if $r_K:= |K|^{1/n}/|B_2^n|^{1/n}$ and
${\cal E}_K = T_K(r_K B_2^n)$ for some volume-preserving $T_K$, then
$\tilde{K}:= T_K^{-1}(K)$ is a linear image of $K$ of the same
volume which is in $M$-position. Assume then that $\tilde{K_1}$ and
$\tilde{K_2}$ are two such images of some bodies $K_1$ and $K_2$ in
${\mathbb R}^n$, and that $K^{\prime}_i$ stands for either
$\tilde{K_i}$ or $(\tilde{K_i})^{\circ}$. Using
(\ref{eq:towardsReverseBM1}) and (\ref{eq:towardsReverseBM2}), we
see that
\begin{align}\label{eq:ReverseBM}
|K^{\prime}_1 + K^{\prime}_2|^{1/n}&\ls c|K^{\prime}_1 +
r_{K^{\prime}_2}B_2^n|^{1/n} \ls c^2|r_{K^{\prime}_1}B_2^n +
r_{K^{\prime}_2}B_2^n|^{1/n}
\\
\nonumber &= c^2\bigl(r_{K^{\prime}_1} +
r_{K^{\prime}_2}\bigr)|B_2^n|^{1/n} = c^2 \bigl(|K^{\prime}_1|^{1/n}
+ |K^{\prime}_2|^{1/n}\bigr).
\end{align}
This means that we have a partial inverse to the Brunn-Minkowski
inequality which holds true for certain affine images of any convex
bodies $K_1, K_2$ and the polars of those images. A direct
consequence of (\ref{eq:ReverseBM}) and Corollary
\ref{corollary-KM2} is the following:
\end{remark}

\begin{corollary}\label{corollary-is-this-correct?}
Let $K$ and $L$ be two convex bodies in ${\mathbb R}^n$ of the same
volume which are in $M$-position. Then,
\begin{equation}N(K,tL)^{1/n}\simeq N(L,tK)^{1/n}\end{equation}
for every $t>0$.
\end{corollary}

\noindent {\it Proof.} Since $tL$ and $tK$ are also in $M$-position for every $t>0$, we have that
\begin{align}
\hspace{0,9cm} N(K,tL)^{1/n}&\simeq \frac{|K + tL|^{\frac{1}{n}}}{|tL|^{\frac{1}{n}}} \simeq
\frac{|K|^{\frac{1}{n}} + t|L|^{\frac{1}{n}}}{t|L|^{\frac{1}{n}}}
\\
\nonumber
&= \frac{t|K|^{\frac{1}{n}} + |L|^{\frac{1}{n}}}{t|K|^{\frac{1}{n}}}\simeq
\frac{|tK + L|^{\frac{1}{n}}}{|tK|^{\frac{1}{n}}}\simeq N(L,tK)^{1/n}.\hspace{1,1cm}\Box
\end{align}
\medskip

Finally, let us remark that, as Pisier notes in \cite{Pisier},
the asymptotic form of the Santal\'{o} inequality and its inverse and
the existence of an $M$-position for any convex body are
interconnected results: if we know that for every centered convex
body $K$ there exists an ellipsoid ${\cal E}_K$ such that
\begin{equation}\max\big\{\log N(K,{\cal E}_K), \log N({\cal E}_K,K), \log N( K^{\circ},{\cal E}_K^{\circ}),
\log N({\cal E}_K^{\circ},K^{\circ})\big\}\ls cn\end{equation} for
some absolute constant $c>0$, then we can prove that
\begin{equation}\label{isomorphic-Santalo}
e^{-2(c+\log\! 8)}s(B_2^n)\ls s(K)\ls e^{2(c+\log\! 8)} s(B_2^n)
\end{equation}
for all centered bodies $K$. Indeed, if ${\cal E}_K$ is an
$M$-ellipsoid for $K$ as above, then from Lemma \ref{lemma-KM2},
\begin{equation*}
\frac{|{\cal E}_K+ K|^{1/n}}{|K|^{1/n}}\ls 2N({\cal E}_K,K)^{1/n}\ls 2e^c \ls
2e^c N(K,{\cal E}_K)^{1/n}\ls 8e^c\frac{|{\cal E}_K+ K|^{1/n}}{|{\cal E}_K|^{1/n}},
\end{equation*}
so $|{\cal E}_K|^{1/n}\ls 8e^c|K|^{1/n}$, and in the
same manner,
\begin{equation}\max\Bigl\{\frac{|K|^{1/n}}{|{\cal E}_K|^{1/n}},\frac{|{\cal E}_K^{\circ}|^{1/n}}{|K^{\circ}|^{1/n}},
\frac{|K^{\circ}|^{1/n}}{|{\cal E}_K^{\circ}|^{1/n}}\Bigr\}\ls 8e^c.\end{equation}
(\ref{isomorphic-Santalo}) now follows.
\bigskip

\medskip

\noindent {\bf Acknowledgements}. The second named author wishes to
thank the A. Sloan Foundation and the US National Science Foundation for
support through the grant DMS-0906150. The third named author is
supported by a scholarship of the University of Athens.

\bigskip

\bigskip

\noindent \textsc{Apostolos Giannopoulos}: Department of
Mathematics, National and Kapodistrian University of Athens,
Panepistimioupolis 157-84, Athens, Greece.

\smallskip

\noindent \textit{E-mail:} \texttt{apgiannop@math.uoa.gr}

\bigskip

\noindent \textsc{Grigoris Paouris}: Department of Mathematics,
Texas A \& M University, College Station, TX 77843 U.S.A.

\smallskip

\noindent \textit{E-mail:} {\tt grigorios.paouris@gmail.com}

\bigskip

\noindent \textsc{Beatrice-Helen Vritsiou}: Department of
Mathematics, National and Kapodistrian University of Athens,
Panepistimioupolis 157-84, Athens, Greece.

\smallskip

\noindent \textit{E-mail:} \texttt{bevritsi@math.uoa.gr}

\end{document}